\documentclass[10pt,twoside]{article}
\usepackage{graphicx}
\usepackage{amsmath}
\usepackage{latexsym}
\usepackage{Latex-document}

\title{\bf  The Complexity of \vskip -2mm
Accurate Floating Point Computation \vskip 6mm}

\author{J. Demmel\vspace*{-0.5cm}\thanks{Mathematics Department and
Computer Science Division, University of California, Berkeley, CA
94720, USA. E-mail: demmel@cs.berkeley.edu}}

\date{\vspace{-8mm}}

\newtheorem{prop}{Proposition}
\newtheorem{conj}{Conjecture}

\newcommand{\bmat}{\left[ \begin{array}}
\newcommand{\emat}{\end{array} \right]}
\newcommand{\size}{{\rm size}}
\newcommand{\diag}{{\rm diag}}

\begin{document}
\maketitle

\thispagestyle{first} \setcounter{page}{697}

\begin{abstract}\vskip 3mm
Our goal is to find accurate and efficient algorithms, when they exist, for
evaluating rational expressions containing floating point numbers,
and for computing matrix factorizations (like LU and the SVD)
of matrices with rational expressions as entries.
More precisely, {\em accuracy} means the relative
error in the output must be less than one (no matter how tiny
the output is), and {\em efficiency} means that the algorithm
runs in polynomial time. Our goal is challenging because
our accuracy demand is much stricter than usual.

The classes of floating point expressions or matrices that we can
accurately and efficiently evaluate or factor depend strongly on
our model of arithmetic:
\begin{enumerate}
\item In the ``Traditional Model'' (TM), the floating point
result of an operation like $a + b$
is $fl(a + b) = (a + b)(1 + \delta)$, where
$|\delta|$ must be tiny.
\item In the ``Long Exponent Model'' (LEM) each floating point number
$x = f \cdot 2^e$ is represented by the pair of integers $(f,e)$,
and there is no bound on the sizes of the exponents $e$ in the input data.
The LEM supports strictly larger classes of expressions or matrices
than the TM.
\item In the ``Short Exponent Model'' (SEM) each floating point number
$x = f \cdot 2^e$ is also represented by $(f,e)$, but
the input exponent sizes are bounded in terms of
the sizes of the input fractions $f$.
We believe the SEM supports strictly more expressions or matrices
than the LEM.
\end{enumerate}

These classes will be
described by factorizability properties of the rational expressions, or
of the minors of the rational matrices. For each such class, we identify
new algorithms that attain our goals of accuracy and efficiency.
These algorithms are often exponentially faster than prior algorithms,
which would simply use a conventional algorithm with sufficiently high precision.

For example, we can factorize
Cauchy matrices, Vandermonde matrices,
totally positive generalized Vandermonde matrices, and
suitably discretized
differential and integral operators
in all three models much more accurately and efficiently than before.
But we provably cannot add $x+y+z$ accurately in the TM,
even though it is easy in the other models.
\vskip 3mm


\end{abstract}

\vskip 12mm

\section{Introduction} \label{section 1}\setzero
\vskip-5mm \hspace{5mm }

We will survey recent progress and describe open problems in the
area of accurate floating point computation, in particular for matrix
computations.  A very short bibliography would include
\cite{demmelkoev99,demmel98,DGESVD,koevthesis,dhillonthesis,barlowdemmel,
demmelkahan,demmelveselic,LAA_Special_Issue}.

We consider the evaluation of multivariate rational functions $r(x)$
of floating point numbers, and matrix computations on rational matrices $A(x)$,
where each entry $A_{ij}(x)$ is such a rational function. Matrix computations
will include computing determinants (and other minors),
linear equation solving, performing Gaussian Elimination (GE) with various kinds
of pivoting, and computing the singular value decomposition (SVD), among others.
Our goals are {\em accuracy} (computing each solution component with tiny
relative error)
and {\em efficiency} (the algorithm should run in time
bounded by a polynomial function of the input size).

We consider three models of arithmetic, defined in the abstract,
and for each one we try to classify rational expressions and matrices
as to whether they can be evaluated or factored accurately and efficiently
(we will say ``compute(d) accurately and efficiently,'' or ``CAE'' for short).

In the
Traditional ``$1 + \delta$'' Model (TM),
we have $fl(a \otimes b) = (a \otimes b)(1 + \delta)$,
$\otimes \in \{ +, -, \times, \div \}$ and $|\delta| \leq \epsilon$, where
$\epsilon \ll 1$ is called {\em machine precision}. It is the conventional model
for floating point error analysis, and means that every floating point result is
computed with a relative error $\delta$ bounded in magnitude by $\epsilon$.
The values of $\delta$ may be arbitrary real (or complex) numbers satisfying
$| \delta | \leq \epsilon$,
so that any algorithm proven to CAE in the TM must work for arbitrary
real (or complex) number inputs and arbitrary real (or complex) $|\delta| \leq \epsilon$.
The size of the input in the TM is the number of floating point words needed
to describe it, independent of $\epsilon$.

The Long Exponent (LEM) and Short Exponent (SEM) models,
which are implementable on a Turing machine, make errors that
may be described by the TM, but their inputs and $\delta$'s are
much more constrained.
Also, we compute the size of the input in the LEM and SEM by counting the number
of bits, so that higher precision and wider range take more bits.

This will mean that problems we can provably CAE
in the TM are a strict subset of those we can CAE in the LEM,
which in turn we conjecture are a strict subset of those we can CAE in the SEM.
In all three models we will describe the classes of rational expressions and rational
matrices in terms of the factorization properties of the expressions, or of the minors
of the matrices.

The reader may wonder why we insist on accurately computing tiny quantities with
small relative error, since in many cases the inputs are themselves uncertain, so
that one could suspect that the inherent uncertainty in the input could make even the
signs of tiny outputs uncertain.
It will turn out that in the TM, the class we can CAE appears to be
identical to the class where all the outputs are in fact accurately determined by
the inputs, in the sense that small relative changes in the inputs cause small
relative changes in the outputs. We make this conjecture more precise in
section 3 below.

There are many ways to formulate the search for efficient and accurate algorithms
\cite{cuckersmale99,smale2,smale3,BCSS,higham96,ReliableComputing,pourelrichards89}. Our approach differs in
several ways.  In contrast to either conventional floating point error analysis \cite{higham96} or the model in
\cite{cuckersmale99}, we ask that even the tiniest results have correct leading digits, and that zero be exact. In
\cite{cuckersmale99} the model of arithmetic allows a tiny absolute error in each operation, whereas in TM we
allow a tiny relative error. Unlike \cite{cuckersmale99} our LEM and SEM are conventional Turing machine models,
with numbers represented as bit strings, and so we can take the cost of arithmetic on very large and very small
numbers (i.e. those with many exponent bits) into precise account. For these reasons we believe our models are
closer to computational practice than the model in \cite{cuckersmale99}. In contrast to \cite{pourelrichards89},
we (mostly) consider the input as given exactly, rather than as a sequence of ever better approximations. Finally,
many of our algorithms could easily be modified to explicitly compute guaranteed interval bounds on the output
\cite{ReliableComputing}.

\section{Factorizability and minors}
\label{section 2}\setzero \vskip-5mm \hspace{5mm }

We show here how to reduce the question of accurate and efficient
matrix computations to accurate and efficient rational expression
evaluation. The connection is elementary, except for the SVD,
which requires an algorithm from \cite{demmelkoev99}.

\begin{prop}
\label{prop_CAEnecessity}
Being able to CAE the absolute value of
the determinant $|\det (A(x))|$ is {\em necessary} to be able to
CAE the following matrix computations on $A(x)$: LU factorization
(with or without pivoting), QR factorization, all the eigenvalues
$\lambda_i$ of $A(x)$, and all the singular values of $A(x)$.
Conversely, being able to CAE {\em all} the minors of $(A(x))$ is
{\em sufficient} to be able to CAE the following matrix
computations on $A(x)$: $A^{-1}$, LU factorization (with or
without pivoting), and the SVD of $A(x)$. This holds in any model
of arithmetic.
\end{prop}
{\bf Proof } First consider necessity. $|\det (A(x))|$ may be
written as the product of diagonal entries of the matrices $L$,
$U$ and $R$ in these factorizations, or as the product of
eigenvalues or singular values. If these entries or values can be
CAE, then so can their product in a straightforward way.

Now consider sufficiency. The statement about $A^{-1}$ is just
Cramer's rule, which only needs $n^2+1$ different minors. The
statement about LU factorization depends on the fact that each
nontrivial entry of $L$ and $U$ is a quotient of minors. The SVD
is more difficult \cite{demmelkoev99}, and depends on the
following two step algorithm: (1) Compute a {\em rank revealing}
decomposition $A = X \cdot D \cdot Y$ where $X$ and $Y$ are
``well-conditioned'' (far from singular in the sense that $\|X\|
\cdot \| X^{-1} \|$ is not too large) and (2) use a
bisection-like algorithm to compute the SVD from $XDY$.

We believe that computing $\det (A(x))$ is actually necessary, not
just $|\det (A(x))|$. The sufficiency proof can be extended to
other matrix computations like the QR decomposition and
pseudoinverse by considering minors of matrices like $\bmat{cc} I
& A \\ A^T & 0 \emat$. Furthermore, if we can CAE the minors of
$C(x) \cdot A(x) \cdot B(x)$, and $C(x)$ and $B(x)$ are
well-conditioned, then we can still CAE a number of matrix
factorizations, like the SVD. The SVD can be applied to get the
eigendecomposition of symmetric matrices, but we know of no
sufficient condition for the accurate and efficient calculation of
eigenvalues of nonsymmetric matrices.

\section{CAE in the traditional model}
\label{section 3}\setzero \vskip-5mm \hspace{5mm}

We begin by giving examples of expressions and matrix computations
that we can CAE in the TM, and then discuss what we cannot do.
The results will depend on details of the axioms we adopt, but
for now we consider the minimal set of operations described in
the abstract.

As long as we only do {\em admissible operations}, namely
multiplication, division, addition of like-signed quantities,
and addition/subtraction of (exact!) input data ($x \pm y$), then
the worst case relative error only grows very slowly, roughly
proportionally to the number of operations. It is when we
subtract two like-signed approximate quantities and significant
cancellation occurs, that the relative error can become large.
So we may ask which problems we can CAE just using only admissible
operations, i.e. which rational expressions factor in such a
way that only admissible operations are needed to evaluate them,
and which matrices have all minors with the same property.

Here are some examples, where we assume that the inputs are arbitrary
real or complex numbers.
(1) The determinant of a Cauchy matrix
$C_{ij} = 1/(x_i + y_j)$ is CAE using the classical expression
$\prod_{i<j} (x_j - x_i)(y_j - y_i)/\prod_{i,j} (x_i + y_j)$,
as is every minor. In fact, changing one line of the classical
GE routine will compute each entry of the LU
decomposition accurately in about the same time as the original
inaccurate version.
(2) We can CAE all minors of sparse matrices,
i.e. those with certain entries fixed at 0 and the rest independent
indeterminates $x_{ij}$, if and only if the undirected bipartite
graph presenting the sparsity structure of the matrix is {\em acyclic};
a one-line change to
GE again renders it accurate. An important special case are bidiagonal
matrices, which arise in the conventional SVD algorithm.
(3) The eigenvalue problem for the second centered difference
approximation to a Sturm-Liouville ODE or elliptic PDE on a
rectangular grid (with arbitrary rectilinear boundaries) can be written
as the SVD of an ``unassembled'' problem $G = D_1 U D_2$ where
$D_1$ and $D_2$ are diagonal (depending on ``masses'' and ``stiffnesses'')
and $U$ is {\em totally unimodular},
i.e. all its minors are $\pm 1$ or 0. Again, a simple change to
GE renders it accurate.

In contrast, one can show that it is impossible in the TM to
add $x+y+z$ accurately in constant time; the proof involves
showing that for {\em any} algorithm the rounding errors $\delta$
and inputs $x,y,z$ can be chosen to have an arbitrarily large
relative error. This depends on the $\delta$'s being permitted
to be arbitrary real numbers in our model.

Vandermonde matrices $V_{ij} = x_i^{j-1}$ are more subtle. Since the product of
a Vandermonde matrix and the Discrete Fourier Transform (DFT)
is Cauchy, and we can compute the SVD of a Cauchy, we can
compute the SVD of a Vandermonde. This fits in our TM model because
the roots of unity in the DFT need only be known approximately,
and so may be computed in the TM model. In contrast, one can
use the result in the last paragraph to show that the inverse of a
Vandermonde cannot be computed accurately. Similarly, polynomial
Vandermonde matrices with $V_{ij} = P_i(x_j)$, $P_i$ a
(normalized) orthogonal polynomial, also permit accurate SVDs,
but probably not inverses.

\section{Adding nonnegativity to the traditional model}
\label{section 4}\setzero \vskip-5mm \hspace{5mm}

If we further restrict the domain of (some) inputs to be nonnegative, then
much more is possible, $x+y+z$ as a trivial example. A more interesting
example are weakly diagonally dominant M-matrices, which arise as discretizations
of PDEs; they must be represented as offdiagonal entries and the row sums.

More interesting is the class of {\em totally positive (TP) matrices}, all of whose
minors are positive. Numerous structure theorems show how to represent such
matrices as products of much simpler TP matrices. Accurate formulas for the
(nonnegative) minors of these simpler matrices combined with the Cauchy-Binet theorem
yield accurate formulas for the minors of the original TP matrix, but typically
at an exponential cost.

An important class of TP matrices where we can do much better are the TP generalized
Vandermonde matrices $G_{ij} = x_i^{\mu_j}$, where the $\mu_j$ form an increasing
nonnegative sequence of integers. $\det (G)$ is known to be the product of
$\prod_{i<j} (x_j-x_i)$ and a {\em Schur function} \cite{macdonald} $s_{\lambda} (x_i)$,
where the sequence $\lambda = (\lambda_j) = (\mu_{n+1-j}-(n-j))$ is called a
{\em partition}.
Schur functions are polynomials with nonnegative integer coefficients, so since their
arguments $x_i$ are nonnegative, they can certainly be computed accurately.
However, straightforward evaluation would have an exponential cost
$O(n^{|\lambda|})$, $|\lambda| = \sum_j \lambda_j$. But by exploiting
combinatorial identities satisfied by Schur functions along with techniques of
divide-and-conquer and memoization, the cost of evaluating the determinant
can be reduced to polynomial time $n^2\prod_j (\lambda_j+1)^2$. The cost of
arbitrary minors and the SVD remains exponential at this time. Note that
the $\lambda_i$ are counted as part of the size of the input in this case.

Here is our conjecture generalizing all the cases we have studied
in the TM. We suppose that $f(x_1,...,x_n)$ is a homogeneous
polynomial, to be evaluated on a domain $\cal D$. We assume that
${\cal D} \subseteq \overline { {\rm int} { \cal D } }$, to avoid
pathological domains. Typical domains could be all tuples of the
real or complex numbers, or the positive orthant.
We say that $f$ satisfies condition $(A)$ (for {\em Accurate}) if
$f$ can be written as a product $f = \prod_m f_m$ where
each factor $f_m$ satisfies
\begin{itemize}
\item $f_m$ is of the form $x_i$, $x_i - x_j$ or $x_i + x_j$, or
\item $|f_m|$ is bounded away from 0 on $\cal D$.
\end{itemize}

\begin{conj}
Let $f$ and $\cal D$ be as above. Then condition $(A)$ is a necessary and
sufficient condition for the existence of an algorithm in the TM model
to compute $f$ accurately on $\cal D$.
\end{conj}

Note that we make no claims that $f$ can be evaluated
efficiently; there are numerous examples where we only know
exponential-time algorithms (doing GE with complete pivoting on a
totally positive generalized Vandermonde matrix).

\section{Extending the TM}
\label{section 5}\setzero \vskip-5mm \hspace{5mm }

So far we have considered the simplest version of the TM, where
(1) we have only the input data, and no additional constants available,
(not even integers, let alone arbitrary rationals or reals),
(2) the input data is given exactly
(as opposed to within a factor of $1+\delta$), and
(3) there is no way to ``round'' a real number to an integer,
and so convert the problem to the LEM or SEM models.
We note that in \cite{cuckersmale99}, (1) integers are available,
(2) the input is rounded, and (3) there is no way to ``round''
to an integer. Changes to these model assumptions will affect
the classes of problems we can solve. For example, if we (quite reasonably)
were to permit exact integers as input, then we could CAE expressions like
$x-1$, and otherwise presumably not. If we went further and permitted
exact rational numbers, then we could also CAE $9x^2-1 = 9(x-\frac{1}{3})(x+\frac{1}{3})$.
Allowing algebraic numbers would make $x^2-2 = (x-\sqrt{2})(x+\sqrt{2})$ CAE.

If inputs were not given exactly, but rather first multiplied by a factor $1+\delta$,
then we could no longer accurately compute
$x \pm y$ where $x$ and $y$ are inputs, eliminating Cauchy matrices and most
others.  But the problems we could solve with exact inputs in the TM still have an
attractive property with inexact inputs:  Small relative changes in the inputs
cause only a small relative change in the outputs, independent of their magnitudes.
The output relative errors may be larger than the input relative error by a factor
called a {\em relative condition number} $\kappa_{rel}$,
which is at most a polynomial function
of $\max (1/ {\rm  rel\_gap} (x_i,\pm x_j))$. Here
${\rm rel\_gap} (x_i, \pm x_j) = |x_i \mp x_j|/( |x_i| + |x_j| )$ is the
{\em relative gap} between inputs $x_i$ and $\pm x_j$, and the maximum is taken
over all expressions $x_i \mp x_j$ where
appearing in $f = \prod_m f_m$.
So if all the input differ in several of their leading digits,
all the leading digits of the outputs are determined accurately.
We note that $\kappa_{rel}$ can be large, depending on $f$ and $\cal D$,
but it can only be unbounded when a relative gap goes to zero.

If a problem has this attractive property, we say that it possesses a
relative perturbation theory.
In practical situations, where only a few leading digits of the inputs $x_i$ are
known, this property justifies the use of algorithms that try to
compute the output as accurately as we do. We state a conjecture very much
like the last one about when a relative perturbation theory exists.

\begin{conj}
Let $f$ and $\cal D$ be as in the last conjecture.
Then condition $(A)$ is a necessary and sufficient condition
for $f$ to have a relative perturbation theory.
\end{conj}

\section{CAE in the long and short exponent models}
\label{section 6}\setzero \vskip-5mm \hspace{5mm }

Now we consider standard Turing machines, where input floating point numbers $x = f \cdot 2^e$ are stored as the
pair of integers $(f,e)$, so the size of $x$ is $\size (x) = \#{\rm bits}(f) + \#{\rm bits}(e)$. We distinguish
two cases, the Long Exponent Model (LEM) where $f$ and $e$ may each be arbitrary integers, and the Short Exponent
Model (SEM), where the length of $e$ is bounded depending on the length of $f$. In the simplest case, when $e=0$
(or lies in a fixed range) then the SEM is equivalent to taking integer inputs, where the complexity of problems
is well understood. This is more generally the case if $\#{\rm bits}(e)$ grows no faster than a polynomial
function of $\#{\rm bits}(f)$.

In particular it is possible to CAE the determinant of an integer
(or SEM) matrix each of whose entries is an independent floating
point number \cite{clarkson}.
This is not possible as far as we know in the LEM,
which accounts for a large complexity gap between the two models.

We start by illustrating some differences between the LEM and SEM,
and then describe the class of problems that we can CAE in the LEM.

First, consider the number of bits in an expression with LEM inputs
can be exponentially larger than the number of bits in the same
expression when evaluated with SEM inputs. For example,
$\size(x \cdot y) \leq \size(x) + \size(y)$ when $x$ and $y$ are
integers, but
$\size(x \cdot y) \leq \size(x) \cdot \size(y)$ when $x$ and $y$ are
LEM numbers:
$(\sum_{i=1}^n 2^{e_i}) \cdot (\sum_{i=1}^n 2^{e'_i})$ has up to
$n^2$ different bit positions to store, each $2^{e_i + e'_j}$, not $2n$.
In other words, LEM arithmetic can encode symbolic algebra, because if
$e_1$ and $e_2$ have no overlapping bits, then we can recover $e_1$ and $e_2$
from the product $2^{e_1} \cdot 2^{e_2} = 2^{e_1 + e_2}$.

Second, the error of many conventional matrix algorithms is typically
proportional to the condition number $\kappa (A) = \|A\| \cdot \|A^{-1}\|$.
This means that a conventional algorithm run with $O(\log \kappa(A))$ extra
bits of precision will compute an accurate answer. It turns out that
if $A(x)$ has rational entries in the SEM model, then
$\log \kappa(A)$ is at most a polynomial function of the input size,
so conventional  algorithms run in high precision will CAE the answer.
However $\log \kappa (A)$ for LEM matrices can be exponentially larger,
so this approach does not work. The simplest example is
$\log \kappa ( \diag (1,2^e) ) = e = 2^{\#{\rm bits}(e)}$.
On the other hand $\log \log \kappa (A(x))$ is a lower bound on the
complexity of any algorithm, because this is a lower bound on the
number of exponent bits in the answer. One can show that
$\log \log \kappa (A(x))$ grows at most polynomially large in the size
of the input.

Finally, we consider the problem of computing an arbitrary bit
in the simple expression $p = \prod_{i=1}^n (1+x_i)$. When the $x_i$
are in the SEM, then $p$ can be computed exactly in polynomial time. However
when the $x_i$ are in the LEM, then one can prove that computing an arbitrary
bit of $p$ is as hard as computing the permanent, a well-known
combinatorially difficult problem.
Here is another apparently simple problem not known to
even be in NP: testing singularity of a floating point matrix. In the
SEM, we can CAE the determinant. But in the LEM, the obvious choice
of a ``witness'' for singularity, a null vector, can have exponentially
many bits in it, even if the matrix is just tridiagonal.
We conjecture that deciding singularity of an LEM matrix is NP-hard.

So how do we compute efficiently in the LEM? The idea is to use
{\em sparse arithmetic}, or to represent only the nonzero bits
in the number. (A long string of 1s can be represented as the
difference of two powers of 2 and similarly compressed). In
contrast, in the SEM one uses {\em dense arithmetic}, storing
all fraction bits of a number. For example, in sparse arithmetic
$2^e+1$ takes $O(\log e)$ bits to store in
sparse arithmetic, but $e$ bits in dense arithmetic.
This idea is exploited in practical floating point computation,
where extra precise numbers are stored as arrays of conventional
floating point numbers, with possibly widely different exponents
\cite{priest}.

Now we describe the class of rational functions that we can
CAE in the LEM.  We say the rational function $r(x)$ is in
factored form if $r(x) = \sum_{i=1}^n p_i (x_1,...,x_k)^{e_i}$,
where each $e_i$ is an integer, and $p_i(x_1,...,x_k)$ is written
as an explicit sum of nonzero monomials.  We say $\size (r)$ is
the number of bits needed to represent it in factored form. Then by
(1) computing each monomial in each $p_i$ exactly,
(2) computing the leading bits of their sum $p_i$ using sparse
arithmetic (the cost is basically sorting the bits), and
(3) computing the leading bits of the product of the $p_i^{e_i}$
by conventional rounded multiplication or division, one can evaluate $r(x)$
accurately in time a polynomial in $\size (r)$ and $\size (x)$.
In other words,
the class of rational expression that we can CAE are those that
we can express in factored form in polynomial space.

Now we consider matrix computations. It follows from the last
paragraph that if each minor $r(x)$ of $A(x)$ can be written
in factored form of a size polynomial in the size of $A(x)$,
then we can CAE all the matrix computations that depend on minors.
So the question is which matrix classes $A(x)$ have all their minors
(or just the ones needed for a particular matrix factorization)
expressible in a factored form no more than polynomially larger than
the size of $A(x)$. The obvious way to write $r(x)$, with the
Laplace expansion, is clearly exponentially larger than $A(x)$,
so it is only specially structured $A(x)$ that will work.

All the matrices that we could CAE in the TM are also possible in the LEM.
The most obvious classes of $A(x)$ that we can CAE in the LEM
that were impossible in the TM are gotten by replacing all the
indeterminates in the TM examples by arbitrary rational expressions
of polynomial size. For example, the entries of an M-matrix can be
polynomial-sized rational expressions in other quantities.
Another class are Green's matrices (inverses of tridiagonals),
which can be thought of as discretized integral operators,
with entries written as $A_{ij} = x_i \cdot y_j$.

The obvious question is whether $A$ each of whose entries is an
independent number in the LEM falls in this class. We conjecture
that it does not, as mentioned before.

\section{Conclusions and open problems}
\label{section 7}\setzero \vskip-5mm \hspace{5mm }

Our goal has been to identify rational expressions (or matrices)
that we can evaluate accurately (or on which we can perform accurate matrix computations),
in polynomial time. Accurately means that we want to get a relative
error less than 1, and polynomial time means in a time bounded by a
polynomial function of the input size.

We have defined three reasonable models of arithmetic, the Traditional Model (TM),
the Long Exponent Model (LEM) and the Short Exponent Model (SEM), and
tried to identify the classes of problems that can or cannot be computed
accurately and efficiently for each model.
The TM can be used as a model to do proofs that also hold in the implementable
LEM and SEM, but since it ignores the structure of floating point numbers
as stored in the computer, it is strictly weaker than either the LEM or SEM.
In other words, there are problems (like adding $x+y+z$) that are provably
impossible in the TM but straightforward in the other two models.

We also believe that the LEM is strictly weaker than the SEM, in the sense
that there appear to be computations (like computing the determinant of
a general, or even tridiagonal, matrix) that are possible in polynomial time
in the SEM but not in the LEM. In the SEM, essentially all problems that
can be written down in polynomial space can be solved in polynomial time.
For the LEM, only expressions that can be written in {\em factored form}
in polynomial space can be computed efficiently in polynomial time.

A number of open problems and conjectures were mentioned in the paper.
We mention just one additional one here: What can be said about
the nonsymmetric eigenvalue problem? In other words, what matrix
properties, perhaps related to minors, guarantee that all eigenvalues of
a nonsymmetric matrix can be computed accurately?


\noindent{\bf Acknowledgements }
  The author acknowledges
Benjamin Diament,
Zlatko Drma\v{c},
Stan Eisenstat,
Ming Gu,
William Kahan,
Ivan Slapni\v{c}ar,
Kresimir Veseli\`{c},
and especially
Plamen Koev
for their collaboration over many years
in developing this material.

\label{lastpage}

\end{document}